\documentclass[12pt]{amsart}
\usepackage{graphicx}
\usepackage{graphics}
\usepackage{amsmath}
\usepackage{amscd}
\usepackage{latexsym}
\usepackage{subfigure}
\usepackage{caption}
\usepackage{hyperref}
\usepackage{placeins}

\begin{document}

\textwidth 5.9in
\textheight 7.9in

\evensidemargin .75in
\oddsidemargin .75in

\newtheorem{Thm}{Theorem}
\newtheorem{Lem}[Thm]{Lemma}
\newtheorem{Cor}[Thm]{Corollary}
\newtheorem{Prop}[Thm]{Proposition}
\newtheorem{Rm}{Remark}

\def\a{{\mathbb a}}
\def\C{{\mathbb C}}
\def\A{{\mathbb A}}
\def\B{{\mathbb B}}
\def\D{{\mathbb D}}
\def\E{{\mathbb E}}
\def\R{{\mathbb R}}
\def\P{{\mathbb P}}
\def\S{{\mathbb S}}
\def\Z{{\mathbb Z}}
\def\O{{\mathbb O}}
\def\H{{\mathbb H}}
\def\V{{\mathbb V}}
\def\Q{{\mathbb Q}}
\def\Cn{${\mathcal C}_n$}
\def\CM{\mathcal M}
\def\CG{\mathcal G}
\def\CH{\mathcal H}
\def\CT{\mathcal T}
\def\CF{\mathcal F}
\def\CA{\mathcal A}
\def\CB{\mathcal B}
\def\CD{\mathcal D}
\def\CP{\mathcal P}
\def\CS{\mathcal S}
\def\CZ{\mathcal Z}
\def\CE{\mathcal E}
\def\CL{\mathcal L}
\def\CV{\mathcal V}
\def\CW{\mathcal W}
\def\IC{\mathbb C}
\def\IF{\mathbb F}
\def\IK{\mathcal K}
\def\IL{\mathcal L}
\def\IP{\bf P}
\def\IR{\mathbb R}
\def\IZ{\mathbb Z}

\title{Corks and exotic ribbons in $B^{4}$}
\author{Selman Akbulut}
\keywords{}

\address{Gökova Geometry Topology Institute,
Akyaka, Mehmet Gokovali Sokak,
No:53, Ula, Muğla, Türkiye}
\email{akbulut.selman@gmail.com}

\subjclass{58D27,  58A05, 57R65}
\date{\today}
\begin{abstract} 
By using corks we construct diffeomorphic ribbon disks $D\subset B^{4}$, which are not smootly isotopic rel boundary to each other. \end{abstract}

\maketitle

\setcounter{section}{-1}

\vspace{-.2in}

\section{Introduction}
Consider the manifold $Q^{4}$, introduced 30 years ago in \cite{a1}, called {\it anticork} in \cite{a2}. $Q$ was obtained from $B^{4}$ by removing a properly imbedded smooth disk $D\subset B^{4}$, bounding the knot $K$ of Figure~\ref{c1}. $Q$ provides an example of an exotic homotopy $S^{1}\times B^{3}$ relative to its boundary. Here we will review the construction of $Q$ by relating it to the cork of \cite{a3}, from which we obtain non-isotopic ribbon disks $D, D'\subset B^{4}$ with same boundary, which are diffeomorphic to each other $(B^{4}, D)\approx (B^{4}, D')$. Then, by using the  infinite order (loose) cork of \cite{a4} and \cite{g}, similarly we produce infinitely many ribbon disks in $B^4$ diffeomorphic but not isotopic to each other relative to boundary. 

\section{Construction }\label{exoticdiff}

\begin{figure}[ht]  \begin{center}
\includegraphics[width=.7\textwidth]{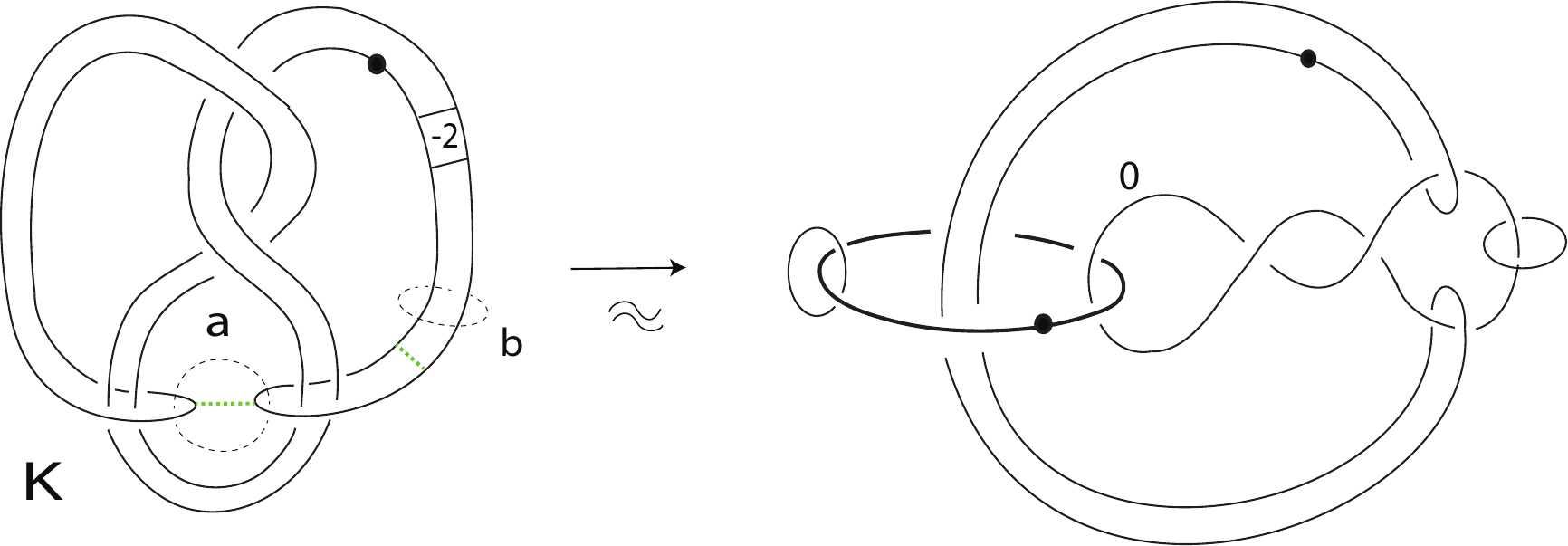}       
\caption{$Q$}      \label{c1} 
\end{center}
 \end{figure}

\begin{figure}[ht]  \begin{center}
\includegraphics[width=.8\textwidth]{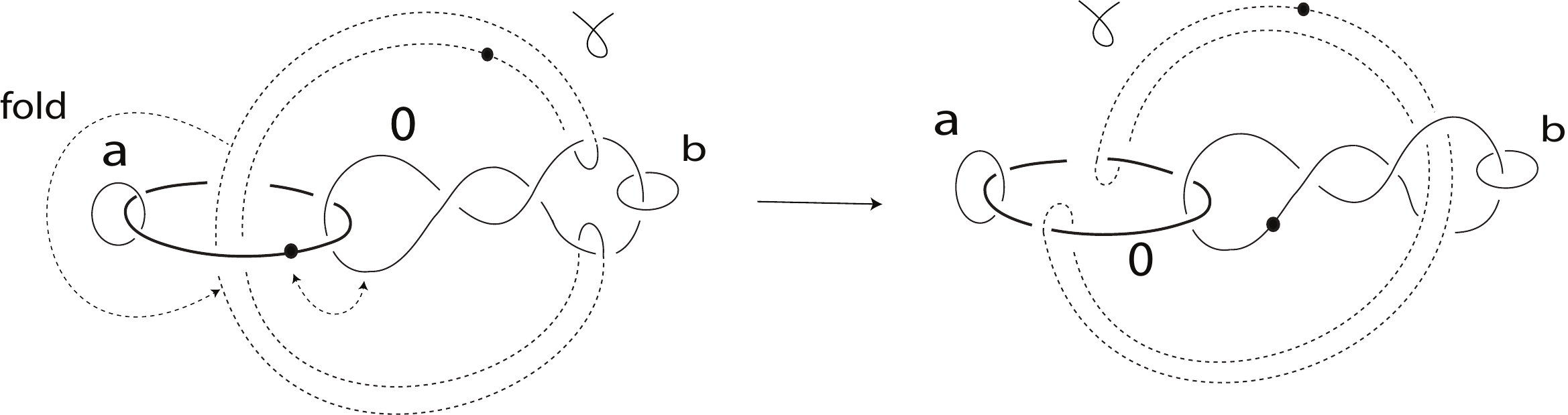}       
\caption{Two pictures of $\gamma \subset \partial B^{4}$}      
\label{c2} 
\end{center}
 \end{figure}

\begin{figure}[ht]  \begin{center}
\includegraphics[width=.88\textwidth]{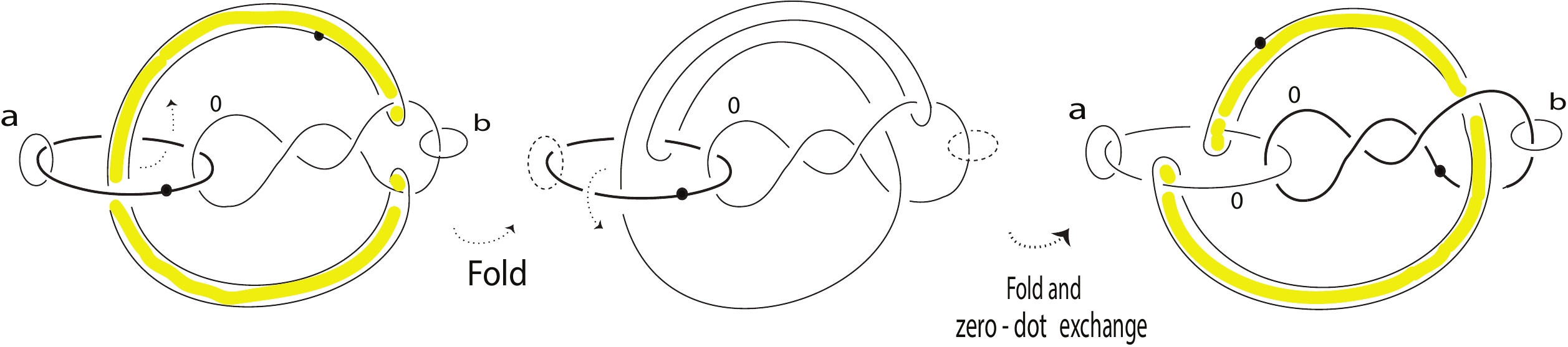}       
\caption{}      \label{c1a} 
\end{center}
 \end{figure}

 Figure~\ref{c1} describes handlebody pictures of $Q$. Recall $Q$ can be obtained in two different ways, by carving $B^{4}$ along either one of two disks $D, D' \subset B^{4}$ bounding the knot $K$ of Figure~\ref{c1}. The two ribbon disks can be seen as follows: Notice that in Figure~\ref{c2} there are two handlebody pictures of $B^{4}$ with the loop $\gamma$ on their boundaries. Check that $\gamma$  corresponds to the knot $K$ in Figure~\ref{c1}. Carving the two $B^{4}$'s along the obvious disks $D$ and $D'$, which $\gamma$ bounds, in either pictures of Figure~\ref{c2}  gives the two descriptions of $Q$.  Notice that even though going from the left to right pictures of Figure~\ref{c2} is just the trivial cork twisting operation of $B^{4}$ (zero-dot exchange), it has a nontrivial affect on the disk $D$, namely it punctures it. $D$ doesn't carry over to the right picture from the left picture, so we have to find a new disk $D'$ which $\gamma$ bounds on the right. To do this we first fold $\gamma$ as shown in Figure~\ref{c2}, so that the clasping part hugs the $2$-handle, then take the obvious disks $D'$ it bounds in the right picture of this Figure~\ref{c2}. Figue~\ref{c1a} shows more detailed picture of this process by shading the disks $D$ and $D'$.

\begin{Thm}
Ribbon disks $D$ and $D'$ can not be smoothly isotopic to each other relative to boundary,  yet there is a diffeomorphism of pairs $(B^{4}, D)\approx (B^{4}, D')$ \end{Thm}

\begin{figure}[ht]  \begin{center}
\includegraphics[width=.78\textwidth]{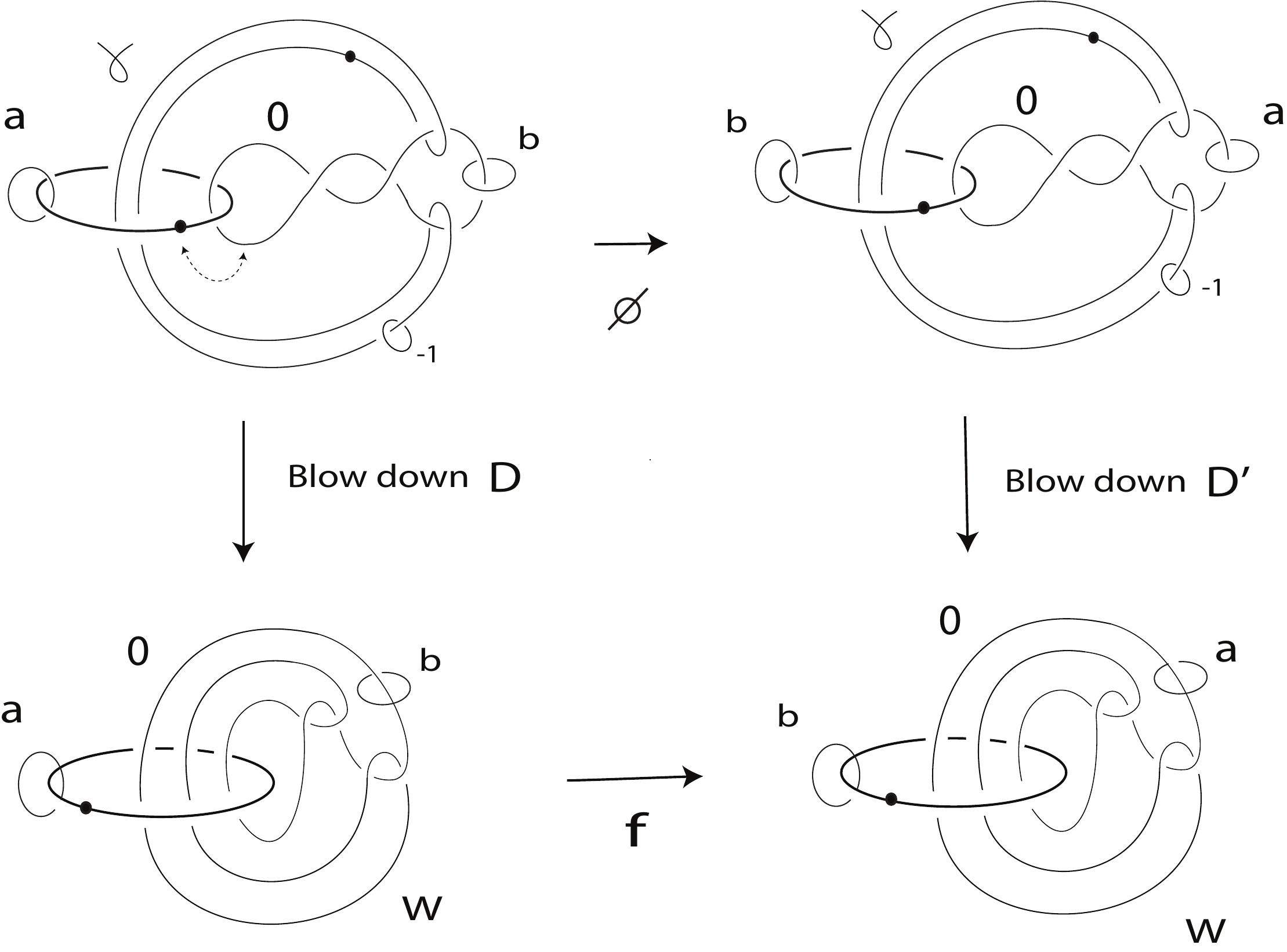}       
\caption{$(W, f)$}      \label{c3} 
\end{center}
 \end{figure}

\proof  It is easy to construct an isotopy $(S^{3},K ) \approx (S^{3},K)$ taking one ribbon arc of $K$ to the other one. This means after initially moving $\partial D =K$ with a concordance near the collar of $\partial B^{4}$, the disks $D$ and $D'$ become isotopic to each other rel bounday. Hence $(B^{4}, D) \approx (B^{4},D')$. On the other hand if the disks $D$ and $D'$ were isotopic to each other rel boundary on the nose, we can blow down $B^{4}$ along $D$ and $D'$ and get a diffeomorphism on the induced contractible manifold $f: W\to W$, as shown in Figure~\ref{c3}, which maps the loop $a$ to $b$ . Notice that $(W,f)$ is the cork appeared in \cite{a3}, and the restriction $f:\partial W \to \partial W$ is its cork twisting map. This is a contradiction $f$ being a cork map.   \qed

 
 \newpage

 \begin{Thm}
There are infinitely many ribbon disks $D_{n}\subset B^{4}$ with the same boundary, which are not smoothly isotopic to each other relative to boundary, yet there are diffeomorphis $(B^{4}, D_{n})\approx (B^{4}, D_{0})$, for all $n=1,2...$
\end{Thm}

The proof of this theorem is similar to the previous theorem, except that here we use the infinite order ``loose'' cork of \cite{a4} rather than the cork $(W,f)$,  discussed in Theorem 1 above.

\newpage

\begin{figure}[ht]
\centering
\begin{minipage}{.35\textwidth}
  \centering
  \includegraphics[width=.8\linewidth]{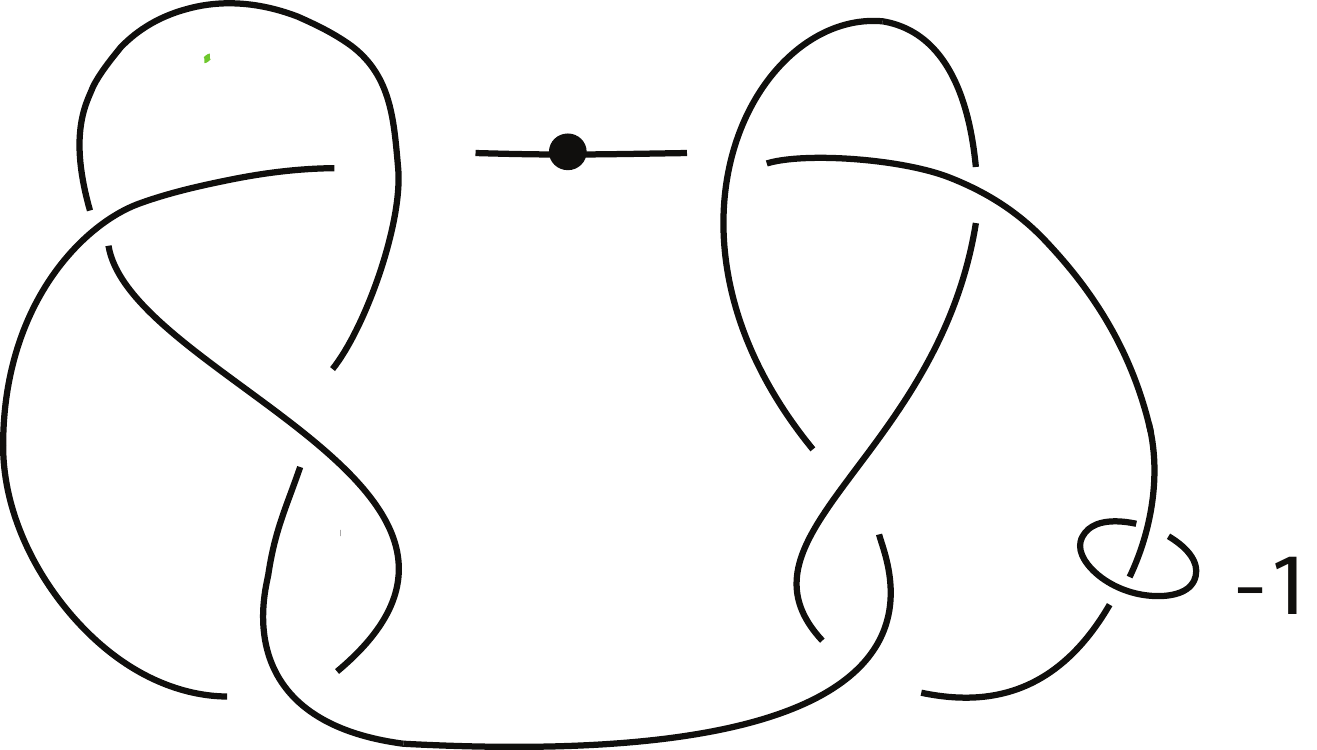}
  \caption{$D$}
  \label{s3}
\end{minipage}%
 \begin{minipage}{.35\textwidth}
  \includegraphics[width=.85\linewidth]{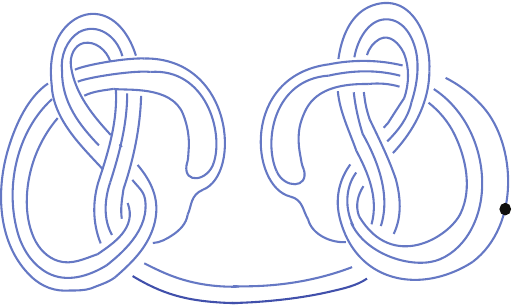}
  \caption{$D_{n}$}
  \label{t6}
\end{minipage}
\end{figure}

\begin{figure}[ht]  \begin{center}
\includegraphics[width=.8\textwidth]{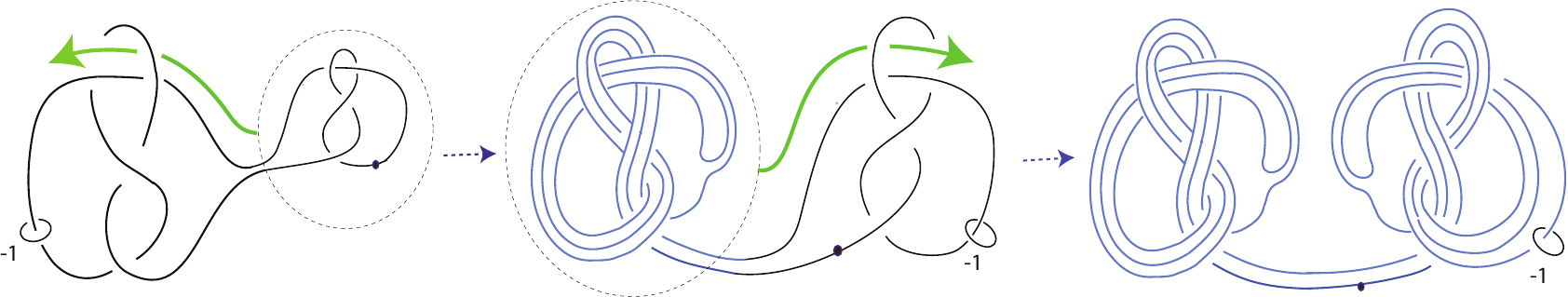}       
\caption{Isotoping the boundary of $D$ to that of $D_{n}$ }      \label{t5} 
\end{center}
 \end{figure}

\begin{figure}[ht]  \begin{center}
\includegraphics[width=.7\textwidth]{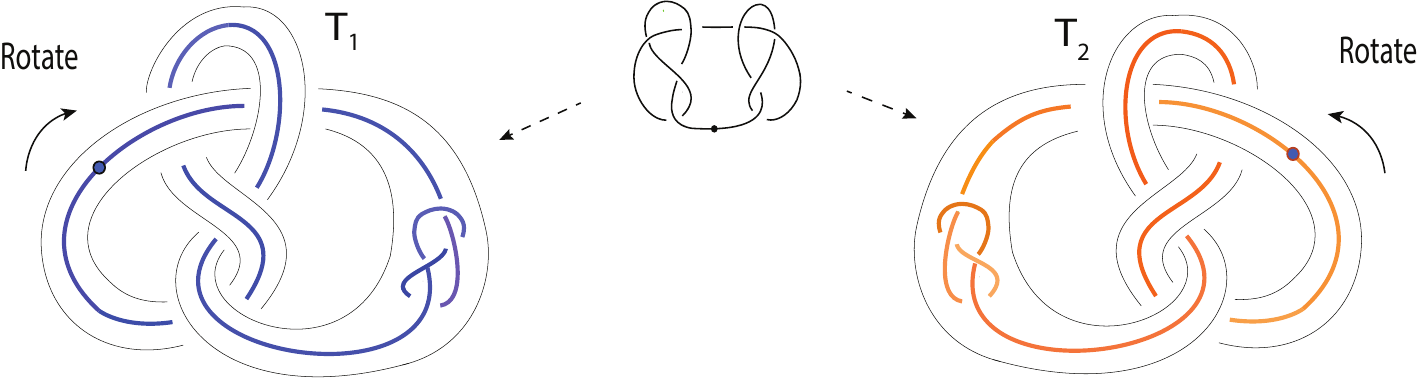}       
\caption{Swallow-follow isotopy}      \label{t7} 
\end{center}
 \end{figure}

\begin{figure}[ht]  \begin{center}
\includegraphics[width=.6\textwidth]{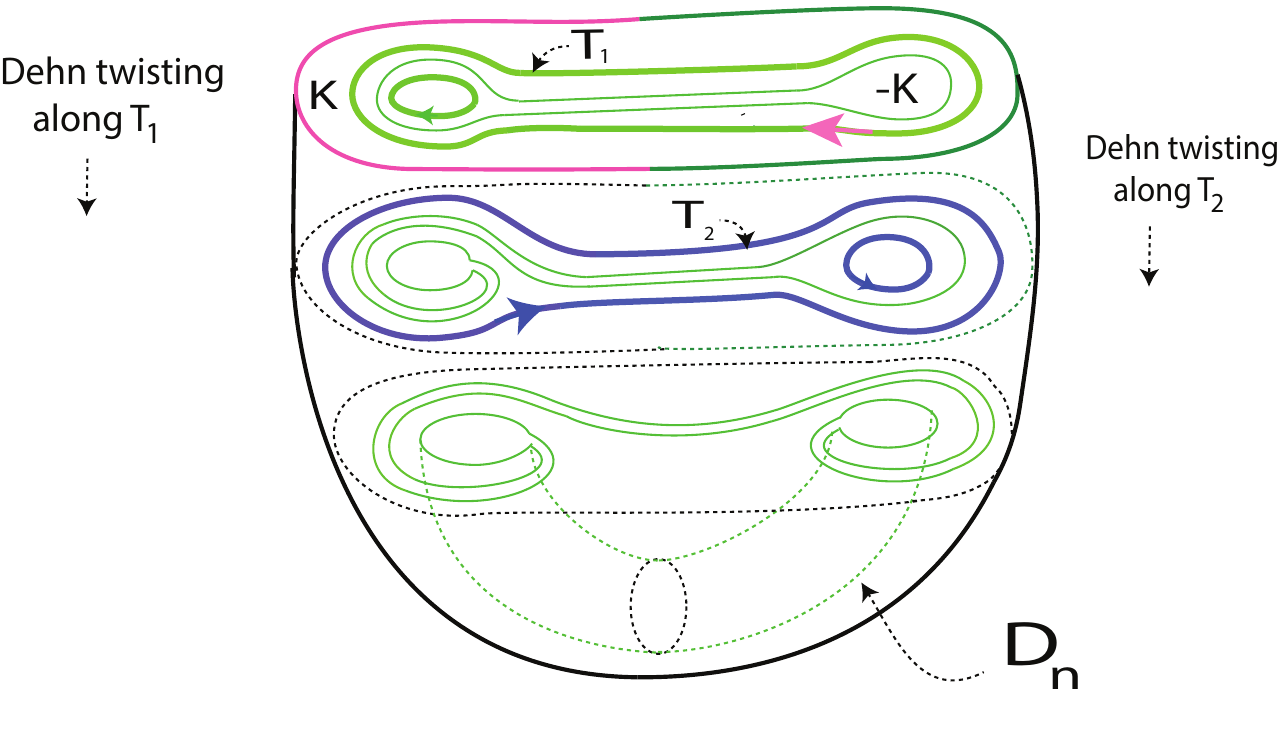}       
\caption{The concordance }      \label{t8} 
\end{center}
 \end{figure}

\proof
Recall the infinite order (loose) cork $(W,f_{n})$ of \cite{a4} and \cite{g}, which is drawn in Figure~\ref{s3}. By \cite{a6} the order $n$ cork-twisting map $f_{n}: \partial W\to \partial W$ is induced from blowing down $B^{4}$, along the ribbon disks $D_{n}\subset B^{4}$ of Figure~\ref{t6}, bounding $K\#-K$, where $K$ is the figure-8 knot and $-K$ is its mirror image (the role of $n$ is self explanatory,  and the picture is drawn for $n=1$ case). So the disks $D_{n}$ can not be isotopic to each other relative to boundary, even though they are diffeomorphic  $(B^{4}, D_{0})\approx (B^{4}, D_{n})$. They become isotopic rel boundary to each other, after isotoping the boundary of one of them in $S^{3}$ as in Figure~\ref{t5} (this is an application of Lemma 2.2 of \cite{jz}). This is so called  ``swallow-follow" isotopy of the knot $K\#-K$, indicated in Figures~\ref{t7}. This isotopy describes a concordance of $K\#-K$, induced by a pair of Dehn twists along the tori $T_{1}$ and $T_{2}$ indicated in Figure~\ref{t8} (see \cite{a6}).  \qed

\end{document}